\documentclass[10pt]{article}
\usepackage{}
\usepackage{indentfirst}
\usepackage{amsfonts}
\usepackage{amsfonts}
\usepackage{amsmath}
\allowdisplaybreaks[4]
\usepackage{bbm}
\usepackage{subfigure}
\usepackage{mathrsfs}
\usepackage{graphicx,amsmath,amsfonts,amssymb,amscd,amsthm}
\usepackage[ urlcolor=blue, colorlinks,linkcolor=black,anchorcolor=blue,citecolor=blue]{hyperref}
\newtheoremstyle{kai}
{10pt} {10pt} {} {} {\bfseries} {.} {.5em} {}
\def\EquationsBySection{\def\theequation
{\thesection.\arabic{equation}}%
\@addtoreset{equation}{section}}

\newcommand\old[1]{}
\newcommand{\pend}{\hfill \thicklines \framebox(6.6,6.6)[l]{}}
\renewenvironment{proof}{\noindent {\it  Proof.} \rm}{\pend}
\newtheorem{theorem}{Theorem}[section]
\newtheorem{lemma}{Lemma}[section]
\newtheorem{corollary}{Corollary}[section]

\newtheorem{remark}{Remark}[section]
\newtheorem{definition}{Definition}[section]

 \EquationsBySection \makeatother


\begin{document}

\pagestyle{plain}
\title
{\bf  Almost sure exponential stabilization of linear stochastic systems driven by a general $\alpha$-stable process via discrete-time state feedback}

\author{\small Xuekang Zhang$^1$\footnote{Corresponding author. E-mail:
xkzhang@ahpu.edu.cn (X. Zhang)}
 \\
\footnotesize{$^1$ School of Mathematics-Physics and Finance, and Key Laboratory of Advanced Perception and Intelligent}\\
\footnotesize{ Control of High-end Equipment, Ministry of Education, Anhui Polytechnic University, Wuhu 241000, China.}
}

\date{}
\maketitle \indent {\bf Abstract.} Given an unstable linear scalar differential equation $dx(t)=\theta x(t)dt$, $\theta>0$, we show that it can be almost surely exponentially stabilized  by a discrete-time feedback control driven by a general $\alpha$-stable process. More precisely, consider the stochastically controlled system
$dx(t)=\theta x(t)dt+bx(\lfloor t/\tau\rfloor\tau)dZ(t)$,
 where $Z(t)$ is a general $\alpha$-stable process with $\alpha\in(1,2)$, and $\lfloor t/\tau\rfloor$ denotes the integer part of $t/\tau$. Unlike the extensively studied driven Brownian-driven case, the $\alpha$-stable process has no finite second moment, so high-order moment methods are no longer applicable. By using fractional moments, characteristic functions, and moment inequalities for $\alpha$-stable stochastic integrals, we prove that the stochastically controlled system is almost surely exponentially stable provided that the inequality
$|b\sigma|^\alpha\left(1-\beta\operatorname{sgn}(b)\tan^2\frac{\pi\alpha}{2}\right)>(\theta+b\mu) \frac{-1}{\Gamma(\alpha)\cos\frac{\alpha\pi}{2}}$ holds, and that $\tau>0$ is sufficiently small.


\section{Introduction}
Feedback plays an essential role in stabilizing dynamical systems. Feedback stabilization can be achieved through both deterministic and stochastic state feedback mechanisms. Unlike deterministic state feedback, which modifies the drift term directly, stochastic state feedback can act through a stochastic channel. To see the difference, consider an unstable scalar linear system
\begin{align*}
\dot{x}(t)=\theta x(t),\quad \theta>0,
\end{align*}
with initial condition $x(0)=x_0$. In the classical deterministic setting, one may introduce a state feedback $u(t)=k x(t)$, leading to
$\dot{x}(t)=\theta x(t)+u(t)=(\theta+k)x(t)$. The controlled system is exponentially stable whenever $k<-\theta$. Hence, deterministic state feedback stabilizes the system by directly modifying the drift term, thereby changing the state trajectory relative to the uncontrolled system. The stabilization problem for deterministic control systems has been extensively investigated (see e.g., Allwright et al. \cite{Allwright2005}, Ebihara et al. \cite{Ebihara2011}, Hagiwara and Araki \cite{Hagiwara1988}). This contrasts with stochastic state feedback, where the control enters through a stochastic channel and affects the noise term rather than the drift term alone. For example,  one may design a linear stochastic state feedback control
$u(t)=bx(t)dB(t)$ based on a scalar Brownian motion $B(t)$, so that the stochastically controlled system becomes
\begin{align}\label{cc}
dx(t)=\theta x(t)dt + b x(t)dB(t),\quad x(0)=x_0,
\end{align}
where $b\neq 0$. It is well known that this system is almost surely exponentially stable whenever  $b^2>2\theta$ (see, e.g., Arnold \cite{Arnold1974}, Has'minskii \cite{Hasminskii1981}, Mao \cite{Mao2007}). Moreover,
\begin{align*}
\mathbb{E}[z(t)]=z(0)e^{\theta t}=x(t),
\end{align*}
which coincides with the solution of the original uncontrolled system. It follows that stochastic feedback achieves sample-path stabilization while preserving the first-moment dynamics of the uncontrolled system. This distinctive property is consistent with the volatility-stabilized market concept in mathematical finance (see, e.g., Fernholz and Karatzas \cite{Fernholz2005}, Shkolnikov \cite{Shkolnikov2013}). Over the past few decades, stabilization and destabilization of stochastic systems driven by Brownian motion  have been studied by many authors (see, e.g., Li et al. \cite{Li2021}, Mao \cite{Mao1994}, Peng and Zhu \cite{Peng2026}, Scheutzow \cite{Scheutzow1993}, Zhang et al. \cite{Zhang2025}).

Continuous-time feedback control is conventionally formulated under the idealized assumption of continuous access to the current state $x(t)$.  In practical digital and networked control systems, however, state measurements are typically available only at discrete time instants. Moreover, continuous observation would be far more expensive and resource-intensive than discrete-time sampling. Therefore, it is more practical to discretize the observation and use the sampled state
$x(\delta_t)$ instead, where $\delta_t=\left\lfloor\frac{t}{\tau}\right\rfloor\tau$,
with $\tau>0$ denoting the duration between two consecutive observations, and $\lfloor \frac{t}{\tau}\rfloor$ is the integer part of $\frac{t}{\tau}$ (see, e.g., ). That is, by choosing a positive constant $\tau$, we need observations only at times $0$, $\tau$, $2\tau$, \ldots. For the above scalar system, this leads to the following stochastically controlled system
\begin{align}\label{modelB}
dx(t)=\theta x(t)dt+bx(\delta_t)dB(t).
\end{align}
Mao \cite{mao2016} used the properties of higher-order moments of normal random variables to show that the controlled system (\ref{modelB}) is almost surely exponentially stable when $b^2>2\theta$ and $\tau$ is sufficiently small. Recently, much effort has been devoted to the discrete-time feedback control for stochastic differential equations (see e.g. Dong \cite{Dong2018}, Li et al. \cite{Li2022}, Liu et al. \cite{Liu2019},  Song et al. \cite{Song2018}, Zhang et al. \cite{Zhang2024}). For instance, Fei et al. \cite{Fei2020} studied the stabilization of highly nonlinear hybrid systems by feedback control based on discrete-time state observations. Guo et al. \cite{Guo2025} stabilized an unstable nonlinear differential equation using discrete-time stochastic feedback control.

Stochastically controlled system noise has typically been assumed to be Gaussian, but such assumptions may be unrealistic in many real-world applications. Many natural phenomena also exhibit heavy-tailed distributions. For instance, Brockmann et al. \cite{Brockmann2006} found that the short-time human travel behavior can be well described by an inverse power law of the distance.  Gabaix et al. \cite{Gabaix2003} showed that stock returns follow a power-law distribution. The generalized central limit theorem implies that, under appropriate conditions, the sum of a sequence of i.i.d. random variables with heavy-tailed distributions will tend to an $\alpha$-stable distribution. Moreover, the Gaussian distribution is a special case of $\alpha$-stable distributions. In general, $\alpha$-stable distributions can be used to model many datasets where some assumptions of the classical central limit theorem, such as finite second moments, are violated, thus offering greater generality. In particular, Mandelbrot \cite{Mandelbrot1963} found that the distributions of monthly wool price changes from 1890 to 1937 follow an
$\alpha$-stable distribution with $\alpha=1.7$. Hence, it is more reasonable to replace the driving Brownian motion with an $\alpha$-stable process for such phenomena.

Motivated by the aforementioned works, in this paper, we consider the following linear stochastically controlled system driven by an $\alpha$-stable L\'evy process:
\begin{align}\label{model}
dx(t) = \theta x(t)\,dt + b x\left(\left\lfloor\frac{t}{\tau}\right\rfloor\tau\right)dL(t),
\end{align}
on $t \geq 0$ with initial value $x(0)=x_0 \in \mathbb{R}$. Here $\{L(t)\}_{t\geq 0}$ is an $\alpha$-stable process with $1<\alpha<2$, defined on the probability space $(\Omega,\mathcal{F},\mathbb{P})$ equipped with a right-continuous and increasing family of $\sigma$-algebras $\{\mathcal{F}_t\}_{t\geq 0}$, and $L(0)=0$. The characteristic function of $L(t)$ is given by
\begin{align*}
\mathbb{E}[e^{iuL(t)}]=\exp\left\{\left[i\mu ut-\sigma^\alpha t|u|^\alpha\left(1-i\beta \operatorname{sgn}(u)\tan\frac{\alpha\pi}{2}\right)\right]\right\},
\end{align*}
where $\mu\in\mathbb{R}$ is the location parameter, $\beta\in[-1,1]$ is the skewness parameter, and $\sigma>0$ is the scale parameter. A natural question is whether the system \eqref{model} can be almost surely exponentially stabilized by the discrete-time feedback.
To the best of our knowledge, this problem has not been
investigated yet. The difficulty in such cases arises from the infinite variance property of $\alpha$-stable process. Consequently, the explicit higher-order moment formulas of normal random variables (up to the 6th order) utilized by Mao \cite{mao2016} cannot be directly applied to the system \eqref{model}. Therefore,  establishing almost sure exponential stabilization of system \eqref{model} is a challenging issue.

\section{Preliminaries}
In this section, we give some preliminaries on $\alpha$-stable distributions, $\alpha$-stable processes, and moment inequalities for $\alpha$-stable stochastic integrals.

 \begin{definition}\label{de2.1}
A random variable $\xi$ is said to follow an $\alpha$-stable distribution,
denoted by $\xi \sim S_\alpha (\sigma,\beta,\mu)$, if it has the characteristic function of the following
form:
 \begin{equation*}
        \Phi_\xi(u)=\mathbb{E}\exp\{iu\xi\}=\left\{ \begin{aligned}
                \exp\left\{-\sigma^\alpha|u|^\alpha\left(1-i\beta \operatorname{sgn}(u)\tan\frac{\alpha\pi}{2}\right)+i\mu u\right\},\quad \textrm{if}\,\,\alpha\neq 1,\\
        \exp\left\{-\sigma |u|\left(1+i\beta\frac{2}{\pi} \operatorname{sgn}(u)\log|u|\right)+i\mu u\right\},\quad \textrm{if}\,\,\alpha= 1,
                                \end{aligned} \right.
        \end{equation*}
where $\alpha\in(0,2]$, $\sigma\in(0,\infty)$, $\beta\in[-1,1]$, and $\mu\in(-\infty, \infty)$ are the index of stability, the scale, skewness, and location parameters, respectively.
\end{definition}

If $\xi \sim S_\alpha (\sigma,\beta,\mu)$ and $\alpha\in(0,2)$, then $\mathbb{E}|\xi|^p<\infty$ when $p\in(0,\alpha)$ and $\mathbb{E}|\xi|^p = \infty$ when $p\in[\alpha,\infty)$. When $\mu=0$, we say $\xi$ is strictly $\alpha$-stable. If in addition $\beta=0$, we call $\xi$ symmetric $\alpha$-stable. Note that $\xi$ is strictly 1-stable($\alpha=1$) if and only if $\beta=0$ (symmetric case). We refer to  Janicki and Weron \cite{Janicki}, and Sato  \cite{Sato} for more details on $\alpha$-stable distributions.

\begin{definition}
An $\mathcal{F}_t$-adapted stochastic process $\{Z_t\}_{t\ge0}$ is called a standard $\alpha$-stable process if \\
(1) $Z_0 = 0$, a.s. \\
(2) $Z_t - Z_s \sim S_\alpha ((t-s)^{\frac{1}{\alpha}}, \beta, 0)$, $t> s \ge 0$; \\
(3) For any finite tume points $0\le s_0 < s_1 < ...<s_m <\infty$, the random variables
$Z_{s_0}, Z_{s_1}-Z_{s_0}, ... , Z_{s_m}-Z_{s_{m-1}}$ are independent.
\end{definition}

Note that $\alpha$-stable process has the self-similarity property $Z_{at}\stackrel{d}{=}a^{\frac{1}{\alpha}}Z_t$ with $a>0$. Every L\'evy motion has a unique modification which is c\`{a}dl\`{a}g (right continuous with left limits) and which is also a L\'evy motion. We will henceforth assume that we are using the c\`{a}dl\`{a}g version of given $\alpha$-stable L\'evy motion. When $\alpha\in(0,1)$, $Z_t$ is a pure jump process with locally finite variation; when $\alpha\in(1,2)$, it is a pure jump process with infinite variation and is a martingale; when $\alpha=1$, it is a pure jump semimartingale. For more path properties of $\alpha$-stable processes, we refer to Sato  \cite{Sato}.

The following lemma provides moment inequalities for $\alpha$-stable stochastic integrals, which will be a crucial tool in the proofs of our main results.

\begin{lemma}\label{le2.2}(see, Lemma 4.4, Long and Qian \cite{Long2013})
Let $f(t)$ be a predictable process satisfying $\int_0^T|f(t)|^\alpha dt<\infty$ almost surely for $T<\infty$, and $F:[0,\infty)\rightarrow[0,\infty)$ be a continuous function. We assume that $(i)$ either $f$ is nonnegative or $Z$ is symmetric, and $(ii)$ there exist positive constants $\lambda_0,C$ and $\alpha_0<\alpha$ such that $F(\lambda\nu)\leq C\lambda^{\alpha_0}F(\nu)$ for all $\nu>0$ and all $\lambda\geq \lambda_0$. Then there exist positive constants $c_1$ and $c_2$ depending only on $\alpha,\alpha_0,\beta,C$, and $\lambda_0$ such that for each $T>0$
\begin{align*}
c_1\mathbb{E}\left[F\left(\left(\int_0^T|f(t)|^\alpha dt\right)^{\frac{1}{\alpha}}\right)\right]\leq \mathbb{E}\left[F\left(\sup\limits_{t\leq T}\left|\int_0^tf(s)dZ_s\right|\right)\right]\leq c_2\mathbb{E}\left[F\left(\left(\int_0^T|f(t)|^\alpha dt\right)^{\frac{1}{\alpha}}\right)\right].
\end{align*}
\end{lemma}
  Another important lemma states that the moments of a random variable can be estimated by relating the fractional moments to the characteristic function.

\begin{lemma}(see, Proposition 6, Fageot et al. \cite{Fageot2017} )\label{le2.1}
For a random variable $X$ with characteristic function $\Phi_X$ and $s\in(0,2)$, we have the relation
\begin{align*}
\mathbb{E}|X|^s=\left(\int_\mathbb{R}\frac{1-\cos u}{|u|^{s+1}}du\right)^{-1}\int_\mathbb{R}\frac{1-\Re(\Phi_X)(u)}{|u|^{s+1}}du\in[0,\infty],
\end{align*}
where $\mathcal{R}(z)$ denotes the real part of $z\in \mathbb{C}$.
\end{lemma}

\section{Almost sure exponential stabilization}
In the section, we study the almost sure exponential stabilization for stochastic system (\ref{model}).

\begin{lemma}\label{lem:2.3}
For any $s\in(1,2)$, we have
\begin{align}
\int_0^1 u^{s-1}\cos u\,du - \sin 1 - (s-1)\cos 1 - (s-1)(s-2)\int_1^\infty u^{s-3}\cos u\,du
&= \Gamma(s)\cos\frac{\pi s}{2}, \label{lem:2.3a}\\[4pt]
\int_0^1 u^{s-1}\sin u\,du + \cos 1 - (s-1)\sin 1 - (s-1)(s-2)\int_1^\infty u^{s-3}\sin u\,du
&= \Gamma(s)\sin\frac{\pi s}{2}. \label{lem:2.3b}
\end{align}
\end{lemma}

\begin{proof} Let $\varepsilon>0$ and $I(\varepsilon)=\int_0^\infty e^{-\varepsilon u}u^{s-1}\cos udu $. Note that
\begin{align}\label{c1.1}
\left|\int_0^1 e^{-\varepsilon u}u^{s-1}\cos udu\right|&\leq \int_0^1 u^{s-1}du=\frac{1}{s}.
\end{align}
Since $\lim\limits_{u\rightarrow\infty}u^{s-1}e^{-\frac{\varepsilon}{2}u}=0$, there exists a constant $\bar{C}>0$ such that $u^{s-1}e^{-\frac{\varepsilon}{2}u}\leq \bar{C}$ for all $u\geq 1$.
Then,
\begin{align}\label{c1.2}
\left|\int_1^\infty e^{-\varepsilon u}u^{s-1}\cos udu\right|&\leq \int_1^\infty e^{-\varepsilon u}u^{s-1}du\nonumber\\
&\leq \bar{C}\int_1^\infty e^{-\frac{\varepsilon}{2} u}\nonumber\\
&=\bar{C}\frac{2}{\varepsilon}e^{-\frac{\varepsilon}{2}}
\end{align}
Combining (\ref{c1.1}) with (\ref{c1.2}) gives that $I(\varepsilon)$ is absolutely convergent. Using $\cos u = \Re(e^{iu})$ and
\begin{align*}
\int_0^\infty e^{-zu}u^{s-1}du=\Gamma(s)z^{-s}, \quad \Re(z)>0.
\end{align*}
we have
\begin{align}\label{c1.3}
I(\varepsilon)=\Re\int_0^\infty e^{-(\varepsilon-i)u}u^{s-1}du=\Re\bigl(\Gamma(s)(\varepsilon-i)^{-s}\bigr).
\end{align}
Let $\varepsilon-i=\sqrt{\varepsilon^2+1}e^{-i \arctan (1/\varepsilon)}$. Then, $(\varepsilon-i)^{-s}=(\varepsilon^2+1)^{-s/2}e^{is\arctan(1/\varepsilon)}$. Since $(\varepsilon^2+1)^{-s/2}\rightarrow 1$ and $\arctan(1/\varepsilon)\rightarrow\frac{\pi}{2}$ as $\varepsilon\rightarrow 0$, it follows that $(\varepsilon-i)^{-s}\rightarrow e^{i\frac{\pi s}{2}}$. Applying  Euler's formula and (\ref{c1.3}), one has
\begin{align}\label{I2.0}
\lim\limits_{\varepsilon\rightarrow 0}I(\varepsilon)=\lim\limits_{\varepsilon\rightarrow 0}\Re\bigl(\Gamma(s)(\varepsilon-i)^{-s}\bigr)=\Re\bigl(\Gamma(s)e^{i\frac{\pi s}{2}}\bigr)=\Gamma(s)\cos\frac{\pi s}{2}.
\end{align}
Note that
\begin{align}\label{I2.1}
I(\varepsilon)&=\int_0^1 e^{-\varepsilon u}u^{s-1}\cos u\,du+\int_1^\infty e^{-\varepsilon u}u^{s-1}\cos u\,du\nonumber\\
&:=I_1(\varepsilon)+I_2(\varepsilon).
\end{align}
For $u\in[0,1]$ and any $\varepsilon>0$, we have $|e^{-\varepsilon u}u^{s-1}\cos u|\le u^{s-1}$, and $\int_0^1 u^{s-1}\,du = \frac{1}{s} < \infty$.
By the dominated convergence theorem,
\begin{align}\label{I2.2}
\lim_{\varepsilon\to0^+} I_1(\varepsilon)
= \int_0^1 \lim_{\varepsilon\to0^+} e^{-\varepsilon u}u^{s-1}\cos u \, du
= \int_0^1 u^{s-1}\cos u \, du .
\end{align}
Now we consider $I_2(\varepsilon)$. Integrating by parts yields that
\begin{align}\label{I2.3}
I_2(\varepsilon) &= e^{-\varepsilon u}u^{s-1}\sin u\biggl|_1^\infty
-\int_1^\infty \sin u\cdot\frac{d}{du}\!\big(e^{-\varepsilon u}u^{s-1}\big)du \nonumber\\
&= -e^{-\varepsilon}\sin 1+\varepsilon\int_1^\infty e^{-\varepsilon u}u^{s-1}\sin u\,du
-(s-1)\int_1^\infty e^{-\varepsilon u}u^{s-2}\sin u\,du .
\end{align}
A further integration by parts gives
\begin{align}\label{I2.4}
\int_1^\infty e^{-\varepsilon u}u^{s-1}\sin udu&=-e^{-\varepsilon u}u^{s-1}\cos u\biggl|_{1}^\infty+\int_1^\infty\cos ud(e^{-\varepsilon u}u^{s-1})\nonumber\\
&=e^{-\varepsilon}\cos 1-\varepsilon\int_1^\infty e^{-\varepsilon u}u^{s-1}\cos udu+(s-1)\int_1^\infty e^{-\varepsilon u}u^{s-2}\cos udu
\end{align}
and
\begin{align}\label{I2.5}
\int_1^\infty e^{-\varepsilon u}u^{s-2}\sin udu
&= -e^{-\varepsilon u}u^{s-2}\cos u\biggl|_1^\infty
+\int_1^\infty \cos u\cdot\frac{d}{du}\!\big(e^{-\varepsilon u}u^{s-2}\big)du \nonumber\\
&= e^{-\varepsilon}\cos 1 - \varepsilon\int_1^\infty e^{-\varepsilon u}u^{s-2}\cos udu + (s-2)\int_1^\infty e^{-\varepsilon u}u^{s-3}\cos udu .
\end{align}
By combining (\ref{I2.3})-(\ref{I2.5}), we can conclude that
\begin{align}\label{I2.6}
I_2(\varepsilon) &= -e^{-\varepsilon}\sin 1 - (s-1)e^{-\varepsilon}\cos 1 + \varepsilon e^{-\varepsilon}\cos 1-\varepsilon^2\int_1^\infty e^{-\varepsilon u}u^{s-1}\cos udu \nonumber\\
&\quad + 2(s-1)\varepsilon\int_1^\infty e^{-\varepsilon u}u^{s-2}\cos udu - (s-1)(s-2)\int_1^\infty e^{-\varepsilon u}u^{s-3}\cos udu.
\end{align}
For $u\in[1,\infty)$ and any $\varepsilon>0$, we have
$|e^{-\varepsilon u}u^{s-3}\cos u| \le u^{s-3}$,
and $\int_1^\infty u^{s-3}\,du = \frac{1}{2-s} < \infty$ (since $s<2$).
The dominated convergence theorem implies that
\begin{align}\label{I2.7}
\lim_{\varepsilon\to0^+} \int_1^\infty e^{-\varepsilon u}u^{s-3}\cos u\,du
= \int_1^\infty \lim_{\varepsilon\to0^+} \bigl(e^{-\varepsilon u}u^{s-3}\cos u\bigr)\,du
= \int_1^\infty u^{s-3}\cos u\,du .
\end{align}
Since $s\in(1,2)$, we have
\begin{align}
\left|\varepsilon^2\int_1^\infty e^{-\varepsilon u}u^{s-1}\cos udu\right|&\leq \varepsilon^2\int_0^\infty e^{-\varepsilon u}u^{s-1}du=\Gamma(s)\varepsilon^{2-s}\label{I2.8},\\
\left|\varepsilon\int_1^\infty e^{-\varepsilon u}u^{s-2}\cos udu\right|&\leq \varepsilon \int_0^\infty e^{-\varepsilon u}u^{s-2}du=\Gamma(s-1)\varepsilon^{2-s}\label{I2.9},
\end{align}
which tend to zero as $\varepsilon\rightarrow 0^+$. By (\ref{I2.6})-(\ref{I2.9}), we get that
\begin{align}\label{I2.10}
\lim_{\varepsilon\to0^+}I_2(\varepsilon)=-\sin 1 - (s-1)\cos 1 - (s-1)(s-2)\int_1^\infty u^{s-3}\cos udu.
\end{align}
It follows from (\ref{I2.1}), (\ref{I2.2}), and (\ref{I2.10}) that
\begin{align}\label{I2.11}
\lim_{\varepsilon\to0^+}I(\varepsilon)=\int_0^1 u^{s-1}\cos u \, du-\sin 1 - (s-1)\cos 1 - (s-1)(s-2)\int_1^\infty u^{s-3}\cos udu.
\end{align}
Combining (\ref{I2.0}) with (\ref{I2.11}) gives that (\ref{lem:2.3a}) holds. Similarly, we can show that (\ref{lem:2.3b}) holds. This completes the proof.
\end{proof}

\begin{lemma} \label{le2.4} For any $p\in(0,\alpha-1)$ and $\alpha\in(1,2)$, let
\begin{align*}
g(\tau)=\frac{2\Gamma(1+p)\sin\frac{\pi p}{2}}{\pi}\int_0^\infty\frac{1-e^{-A_1(\tau)u^\alpha}\cos\left(A_2(\tau)u+A_3(\tau)u^\alpha\right)}{u^{p+1}}du,
\end{align*}
where $A_1(\tau)=|b\sigma|^\alpha\frac{e^{\alpha \theta\tau}-1}{\theta\alpha}$, $A_2(\tau)=e^{\theta\tau}+\frac{b\mu}{\theta}(e^{\theta\tau}-1)$, $A_3(\tau)=\widetilde{\beta}A_1(\tau)$, and $\widetilde{\beta}=\beta\operatorname{sgn}(b)\tan\frac{\pi\alpha}{2}$. Then,
 we have\\
 (i) $g(0)=1$, $\lim\limits_{\tau\rightarrow \infty}g(\tau)=\infty$;\\
 (ii)  $g(\tau)$ is continuous on $[0,\infty)$ and is differentiable on $(0,\infty)$;\\
  (iii)
\begin{align}\label{g0}
g'(0^+)&=\frac{2\Gamma(1+p)\sin\frac{\pi p}{2}}{\pi}\left(|b\sigma|^\alpha\Gamma(\alpha-p)\cos\frac{(\alpha-p)\pi}{2}-\widetilde{\beta}|b\sigma|^\alpha\Gamma(\alpha-p)\sin\frac{(\alpha-p)\pi}{2}\right.\nonumber\\
&\qquad\left.-(\theta+b\mu)\Gamma(1-p)\cos\frac{p\pi}{2}\right).
\end{align}
\end{lemma}

\begin{proof}
Note that
\begin{align}\label{G}
\int_0^\infty\frac{1-\cos u}{u^{p+1}}du=\frac{\pi}{2\Gamma(p+1)\sin\frac{\pi p}{2}}.
\end{align}
 A direct computation shows that $g(0)=1$. Let $\tau_0\geq 0$ be arbitrary and sequence $\{\tau_n\}_{n=1}^\infty\subset[0,\infty)$ with $\lim\limits_{n\rightarrow\infty}\tau_n=\tau_0$. Then there exists $T>0$ such that $\tau_n,\tau_0\in[0,T]$ for all $n$. Let \begin{align*}
 f(\tau,u)=\frac{1-e^{-A_1(\tau)u^\alpha}\cos\left(A_2(\tau)u+A_3(\tau)u^\alpha\right)}{u^{p+1}}.
 \end{align*}
  Clearly, $f(\tau,u)\geq 0$ and $\lim_{n\rightarrow\infty}f(\tau_n,u)=f(\tau_0,u)$. Using the elementary inequalities
\begin{align*}
1-\cos u\leq \frac{u^2}{2},\quad 1-e^{-u}\leq u, \quad u\geq 0,
\end{align*}
one can verify that for all $\tau_n\in[0,T]$,
\begin{align*}
f(\tau_n,u)\leq h(u)
\end{align*}
where
\[ h(u)=
\begin{cases}
\frac{|b\sigma|^\alpha(e^{\alpha\theta T}-1)}{\alpha\theta}u^{\alpha-p-1}+\frac{1}{2}\left(e^{\theta T}+\frac{|b\mu|(e^{\theta T}-1)}{\theta}+|\widetilde{\beta}|\frac{|b\sigma|^\alpha(e^{\alpha\theta T}-1)}{\alpha\theta}\right)^2u^{1-p}, & 0 < u \leq 1, \\
\frac{2}{u^{p+1}}, & u \geq 1.
\end{cases}
\]
Obviously,
\begin{align*}
\int_0^\infty h(u)du=\frac{|b\sigma|^\alpha(e^{\alpha\theta T}-1)}{\alpha\theta(\alpha-p)}+\frac{\left(e^{\theta T}+\frac{|b\mu|(e^{\theta T}-1)}{\theta}+|\widetilde{\beta}|\frac{|b\sigma|^\alpha(e^{\alpha\theta T}-1)}{\alpha\theta}\right)^2}{2(2-p)}+\frac{2}{p}.
\end{align*}
By Lebesgue's dominated convergence theorem, it follows that $\lim_{\tau\rightarrow\tau_0}g(\tau)=g(\tau_0)$. Since $\tau_0$ is arbitrary, $g(\tau)$ is continuous on $[0,\infty)$. For any sequence $\{\nu^{(1)}_n\}_{n=1}^\infty\subset[0,\infty)$ with $\lim\limits_{n\rightarrow\infty}\nu^{(1)}_n=\infty$, by Fatou's lemma,
\begin{align*}
\liminf\limits_{n\rightarrow\infty}g(\nu^{(1)}_n)&=\frac{2\Gamma(1+p)\sin\frac{\pi p}{2}}{\pi}\liminf\limits_{n\rightarrow\infty}\int_0^\infty f(\nu^{(1)}_n,u)du\\
&\geq \frac{2\Gamma(1+p)\sin\frac{\pi p}{2}}{\pi}\int_0^\infty \liminf\limits_{n\rightarrow\infty}f(\nu^{(1)}_n,u)du\\
&=\frac{2\Gamma(1+p)\sin\frac{\pi p}{2}}{\pi}\int_0^\infty\frac{1}{u^{p+1}}=\infty.
\end{align*}
This implies that $\lim\limits_{\tau\rightarrow\infty}g(\tau)=\infty$.

From the definition of $f(\tau,u)$, we have
\begin{align*}
\partial_\tau f(\tau,u) &= \frac{e^{-A_1(\tau) u^\alpha} }{u^{p+1}}\left[|b\sigma|^\alpha e^{\alpha\theta\tau} u^\alpha \cos\left(A_2(\tau)u+A_3(\tau)u^\alpha\right)\right.\\
&\quad\left.-\left( (\theta +b\mu)e^{\theta\tau} u+ |b\sigma|^\alpha e^{\alpha\theta\tau} \widetilde{\beta}u^\alpha\right) \sin\left(A_2(\tau)u+A_3(\tau)u^\alpha\right) \right].
\end{align*}
Fix an arbitrary $\tau_1>0$. For every $\tau\in[\frac{\tau_1}{2}, \frac{3\tau_1}{2}]$, we have
\begin{align*}
|\partial_\tau f(\tau,u)| \leq H(u),
\end{align*}
where  $\delta=|b\sigma|^\alpha \frac{e^{\alpha\theta\tau_1/2} - 1}{\alpha\theta}$, and
\begin{align*}
H(u)=|b\sigma|^\alpha e^{\frac{3}{2}\alpha\theta\tau_1}\left(1+\left|\widetilde{\beta}\right|\right)u^{\alpha-1-p} e^{-\delta u^\alpha}+ |\theta+b\mu| e^{\frac{3}{2}\theta\tau_1}u^{-p}e^{-\delta u^\alpha}.
\end{align*}
Indeed,
\begin{align*}
\int_0^\infty H(u)du=|b\sigma|^\alpha e^{\frac{3}{2}\alpha\theta\tau_1}\left(1+\left|\widetilde{\beta}\right|\right)\frac{\delta^{-\frac{\alpha-p}{\alpha}}}{\alpha}\Gamma\left(\frac{\alpha-p}{\alpha}\right)+|\theta+b\mu| e^{\frac{3}{2}\theta\tau_1}\frac{\delta^{-\frac{1-p}{\alpha}}}{\alpha}\Gamma\left(\frac{1-p}{\alpha}\right).
\end{align*}
Hence, we can conclude that  $g(\tau)$ is differentiable at $\tau_1$. Since $\tau_1>0$ is arbitrary, $g(\tau)$ is differential on $(0,\infty)$ and for all $\tau>0$,
 \begin{align}\label{u0}
g'(\tau)&=\frac{2\Gamma(1+p)\sin\frac{\pi p}{2}}{\pi}\left[|b\sigma|^\alpha e^{\alpha\theta\tau}\int_0^\infty e^{-A_1(\tau)u^\alpha}\cos\left(A_2(\tau)u+A_3(\tau)u^\alpha\right)u^{\alpha-p-1}du\right.\nonumber\\
&\left.\quad-\widetilde{\beta}|b\sigma|^\alpha e^{\alpha\theta\tau}\int_0^\infty e^{-A_1(\tau)u^\alpha}\sin\left(A_2(\tau)u+A_3(\tau)u^\alpha\right)u^{\alpha-p-1}du\right]\nonumber\\
&\left.\quad-(\theta+b\mu) e^{\theta\tau}\int_0^\infty e^{-A_1(\tau)u^\alpha}\sin\left(A_2(\tau)u+A_3(\tau)u^\alpha\right)u^{-p}du\right.\nonumber\\
 &=:\frac{2\Gamma(1+p)\sin\frac{\pi p}{2}}{\pi}\left[|b\sigma|^\alpha e^{\alpha\theta\tau}J_1(\tau)-\widetilde{\beta}|b\sigma|^\alpha e^{\alpha\theta\tau}J_2(\tau)-(\theta+b\mu) e^{\theta\tau}J_3(\tau)\right].
\end{align}
Next, we study the asymptotic behavior of $J_i(\tau)$, $i=1,2,3$. For $J_1(\tau)$, we have
\begin{align}\label{J1}
J_1(\tau)&=\int_0^1 e^{-A_1(\tau)u^\alpha}u^{\alpha-p-1}\cos\left(A_2(\tau)u+A_3(\tau)u^\alpha\right)du+\int_1^\infty e^{-A_1(\tau)u^\alpha}u^{\alpha-p-1}\cos\left(A_2(\tau)u+A_3(\tau)u^\alpha\right)du\nonumber\\
&=:J_{11}(\tau)+J_{12}(\tau).
\end{align}
Denote $G(u,\tau)=e^{-A_1(\tau)u^\alpha}u^{\alpha-p-1}$, $u\geq 0$. Observe that
\begin{align}\label{u1}
|G(u,\tau)\cos\left(A_2(\tau)u+A_3(\tau)u^\alpha\right)|\leq u^{\alpha-p-1},\quad u\in[0, 1).
\end{align}
Since $\int_0^1 u^{\alpha-p-1}=\frac{1}{\alpha-p}$, by (\ref{u1}) and dominated convergence theorem, it follows that
\begin{align}\label{J11}
\lim\limits_{\tau\rightarrow 0^+}J_{11}(\tau)=\int_0^1 \lim\limits_{\tau\rightarrow 0^+} e^{-A_1(\tau)u^\alpha}u^{\alpha-p-1}\cos\left(A_2(\tau)u+A_3(\tau)u^\alpha\right)du=\int_0^1u^{\alpha-p-1}\cos udu.
\end{align}
Next, we show that
\begin{align}\label{u6}
\lim\limits_{\tau\rightarrow 0^+}J_{12}(\tau)=-\sin1-(\alpha-p-1)\cos 1-(\alpha-1-p)(\alpha-p-2)\int_1^\infty u^{\alpha-p-3}\cos udu.
\end{align}
We consider two cases according to the sign of $\widetilde{\beta}$.

Case (i): $\widetilde{\beta}\geq 0$.  By integrate by parts, we obtain
\begin{align}\label{u2}
J_{12}(\tau)&=\int_1^\infty \frac{G(u,\tau)}{D(u,\tau)}d\sin\left(A_2(\tau)u+A_3(\tau)u^\alpha\right)\nonumber\\
&=-\frac{G(1,\tau)}{D(1,\tau)}\sin \left(A_2(\tau)+A_3(\tau)\right)-\int_1^\infty\sin\left(A_2(\tau)u+A_3(\tau)u^\alpha\right) \frac{\partial }{\partial  u}\left(\frac{G(u,\tau)}{D(u,\tau)}\right)du\nonumber\\
&=-\frac{G(1,\tau)}{D(1,\tau)}\sin \left(A_2(\tau)+A_3(\tau)\right)-\frac{\cos\left(A_2(\tau)+A_3(\tau)\right)}{D(1,\tau)}\left( \frac{D(1,\tau)G_u(1,\tau)-G(1,\tau)D_u(1,\tau)}{(D(1,\tau))^2} \right)\nonumber\\
&\quad-\int_1^\infty \cos\left(A_2(\tau)u+A_3(\tau)u^\alpha\right)Q(u,\tau)du\nonumber\\
&=:-J_{12}^{(1)}(\tau)-J_{12}^{(2)}(\tau)-J_{12}^{(3)}(\tau),
\end{align}
where
\begin{align*}
D(u,\tau)&=A_2(\tau)+\alpha A_3(\tau)u^{\alpha-1},\\
D_u(u,\tau)&=\alpha(\alpha-1)A_3(\tau)u^{\alpha-2},\\
D_{uu}(u,\tau)&=\alpha(\alpha-1)(\alpha-2)A_3(\tau)u^{\alpha-3},\\
G(u,\tau)&=e^{-A_1(\tau)u^\alpha }u^{\alpha-p-1},\\
G_u(u,\tau)&=e^{-A_1(\tau)u^\alpha }u^{\alpha-p-2}(\alpha-p-1-\alpha A_1(\tau)u^\alpha),\\
G_{uu}(u,\tau)&=e^{-A_1(\tau)u^\alpha }u^{\alpha-p-3}\left[(\alpha-p-1)(\alpha-p-2)-\alpha A_1(\tau)u^\alpha (3\alpha-2p-3)+\alpha^2(A_1(\tau))^2u^{2\alpha}\right],\\
Q(u,\tau)&=\frac{G_{uu}(u,\tau)}{(D(u,\tau))^2}-\frac{3G_u(u,\tau)D_u(u,\tau)+G(u,\tau)D_{uu}(u,\tau)}{(D(u,\tau))^3}+\frac{3G(u,\tau)(D_u(u,\tau))^2}{(D(u,\tau))^4}.
\end{align*}
Obvious that
\begin{align}\label{u3}
\lim\limits_{\tau\rightarrow 0^+}J_{12}^{(1)}(\tau)=\sin1,\quad
\lim\limits_{\tau\rightarrow 0^+}J_{12}^{(2)}(\tau)=(\alpha-p-1)\cos 1,
\end{align}
Since $\lim_{\tau\rightarrow 0^+}A_2(\tau)=1$, there exists a positive constant $\widetilde{\tau}$ such that for all $0<\tau<\widetilde{\tau}$ we have $|A_2(\tau)-1|< \frac{1}{2}$. This implies that
\begin{align}\label{A2}
A_2(\tau)> \frac{1}{2}, \quad \tau\in(0,\widetilde{\tau}).
 \end{align}
Because $\widetilde{\beta}\geq 0$, we have $A_3(\tau)\geq 0$, and therefore for all $\tau\in(0,\widetilde{\tau})$ and $u\geq 1$,
\begin{align*}
D(u,\tau)&\geq A_2(\tau)> \frac{1}{2}.
\end{align*}
Hence, for any $\tau\in(0,\widetilde{\tau})$ and $u\geq 1$, using $y^me^{-y}\leq\left(\frac{m}{e}\right)^{m}$, $m=0,1,2$, we get
\begin{align}\label{u4}
\left|\cos\left(A_2(\tau)u+A_3(\tau)u^\alpha\right)Q(u,\tau)\right|&\leq 4\left[|(\alpha-p-1)(\alpha-p-2)|+\frac{\alpha|3\alpha-2p-3|}{e}+\frac{4\alpha^2}{e^2}  \right]u^{\alpha-p-3}\nonumber\\
&\quad+24\alpha(\alpha-1)|\widetilde{\beta}|\left(\frac{|\alpha-p-1|}{e}+\frac{4\alpha}{e^2}\right)u^{\alpha-p-3}\nonumber\\
&\quad +8\alpha(\alpha-1)|\alpha-2||\widetilde{\beta}|\frac{1}{e}u^{\alpha-p-4}+48\alpha^2(\alpha-1)^2\widetilde{\beta}^2\frac{1}{e^2}u^{\alpha-p-5}\nonumber\\
&\leq K u^{\alpha-p-3},
\end{align}
where
\begin{align*}
K&=4\left[|(\alpha-p-1)(\alpha-p-2)|+\frac{\alpha|3\alpha-2p-3|}{e}+\frac{4\alpha^2}{e^2}  \right]+24\alpha(\alpha-1)|\widetilde{\beta}|\left(\frac{|\alpha-p-1|}{e}+\frac{4\alpha}{e^2}\right)\nonumber\\
&\quad +8\alpha(\alpha-1)|\alpha-2||\widetilde{\beta}|\frac{1}{e}+48\alpha^2(\alpha-1)^2\widetilde{\beta}^2\frac{1}{e^2}.
\end{align*}
Since $\int_1^\infty Ku^{\alpha-p-3}=\frac{K}{2-\alpha+p}$, applying (\ref{u4}) together with the dominated convergence theorem yields that
\begin{align}\label{u5}
\lim\limits_{\tau\rightarrow 0^+}\int_1^\infty \cos\left(A_2(\tau)u+A_3(\tau)u^\alpha\right)Q(u,\tau)du&=(\alpha-1-p)(\alpha-p-2)\int_1^\infty u^{\alpha-p-3}\cos udu.
\end{align}
From (\ref{u2}), (\ref{u3}) and (\ref{u5}), we can conclude that (\ref{u6}) holds.

Case (ii):  $\widetilde{\beta}<0$. Since $\lim_{\tau\rightarrow 0^+}A_3(\tau)=0$, there exists a positive constant $\widetilde{\tau}_1$ such that for all $0<\tau<\widetilde{\tau}_1$ we have $|A_3(\tau)|< \frac{1}{2\alpha}$. Set $\widehat{\tau}=\min\{\widetilde{\tau},\widetilde{\tau}_1 \}$. Then, for all $\tau\in(0,\widehat{\tau})$, by (\ref{A2}), we have $D(1,\tau)>0$.
Clearly, $D(u,\tau)$ is strictly decreasing with respect to $u$ on $[1,\infty)$ and $\lim_{u\rightarrow\infty}D(u,\tau)=-\infty$. By the Intermediate Value Theorem, there exists a unique $u_0(\tau)>1$ such that $D(u_0,\tau)=0$, namely,
\begin{align*}
u_0(\tau)=\left(-\frac{A_2(\tau)}{\alpha\widetilde{\beta}}\right)^{\frac{1}{\alpha-1}}\left(A_1(\tau)\right)^{-\frac{1}{\alpha-1}}.
\end{align*}
Since \(\widetilde{\beta}<0\), one can easily verify that \(\lim_{\tau\to0^+} u_0(\tau)=\infty\). Consequently, there exists $\widetilde{\tau}_2\in(0,\widehat{\tau})$ such that $u_0(\tau)\geq 2$ for all $\tau\in(0,\widetilde{\tau}_2)$. Define $M_\tau=\frac{u_0(\tau)}{2}$, $\tau\in(0,\widetilde{\tau}_2)$.
Note that for all $\tau\in(0,\widetilde{\tau}_2)$, we have
\begin{align}\label{J12}
J_{12}(\tau)&=\int_1^{M_\tau} e^{-A_1(\tau)u^\alpha}u^{\alpha-p-1}\cos\left(A_2(\tau)u+A_3(\tau)u^\alpha\right)du+\int_{M_\tau}^\infty e^{-A_1(\tau)u^\alpha}u^{\alpha-p-1}\cos\left(A_2(\tau)u+A_3(\tau)u^\alpha\right)du\nonumber\\
&=:\widetilde{J}_1(\tau)+\widetilde{J}_2(\tau).
\end{align}
Arguing exactly as in the proof of (\ref{u2}), we get that for all $\tau\in(0,\widetilde{\tau}_2)$,
\begin{align}\label{u7}
\widetilde{J}_1(\tau)&=-J_{12}^{(1)}(\tau)-J_{12}^{(2)}(\tau)+\widetilde{J}_1^{(1)}(\tau)+\widetilde{J}_1^{(2)}(\tau)-\widetilde{J}_1^{(3)}(\tau),
\end{align}
where
\begin{align*}
\widetilde{J}_1^{(1)}(\tau)&=\frac{\cos\left(A_2(\tau)M_\tau+A_3(\tau)M_\tau^\alpha\right)}{D(M_\tau,\tau)}\left( \frac{D(M_\tau,\tau)G_u(M_\tau,\tau)-G(M_\tau,\tau)D_u(M_\tau,\tau)}{(D(M_\tau,\tau))^2} \right),\nonumber\\
\widetilde{J}_1^{(2)}(\tau)&=\frac{G(M_\tau,\tau)}{D(M_\tau,\tau)}\sin\left( A_2(\tau)M_\tau+A_3(\tau)M_\tau^\alpha\right),\nonumber\\
\widetilde{J}_1^{(3)}(\tau)&=\int_1^{M_\tau} \cos\left(A_2(\tau)u+A_3(\tau)u^\alpha\right)Q(u,\tau)du.
\end{align*}
Applying the elementary inequality  $\lim_{y\rightarrow\infty}y^{k_1}e^{-K_2y}=0$ (for any $k_1,k_2>0$), we get that for all $\tau\in(0,\widetilde{\tau}_2)$,
\begin{align*}
|G(M_\tau,\tau)|=\left(2^{-1}(A_1(\tau))^{-\frac{1}{\alpha-1}}\left(-\frac{A_2(\tau)}{\alpha\widetilde{\beta}}\right)^{\frac{1}{\alpha-1}}\right)^{\alpha-p-1}e^{-2^{-\alpha}\left(-\frac{A_2(\tau)}{\alpha\widetilde{\beta}}\right)^{\frac{\alpha}{\alpha-1}}(A_1(\tau))^{-\frac{1}{\alpha-1}}}\rightarrow 0,
\end{align*}
as $\tau\rightarrow0^+$. Similarly, we have
\begin{align*}
|G_u(M_\tau,\tau)|\rightarrow 0,\quad |G(M_\tau,\tau)D_u(M_\tau,\tau)|\rightarrow 0\quad \mathrm{as} \quad \tau\rightarrow 0^+.
\end{align*}
Because $D(u,\tau)$ is strictly decreasing with respect to $u$ on $[1,\infty)$, we obtain that
\begin{align*}
D(M_\tau,\tau)&\geq \frac{1}{2}\left(1-\frac{1}{2^{\alpha-1}}\right),\quad \tau\in(0,\widetilde{\tau}_2).
\end{align*}
Hence, we immediately obtain that
\begin{align}\label{u8}
\lim\limits_{\tau\rightarrow 0^+}\widetilde{J}_1^{(1)}(\tau)=0,\quad \lim\limits_{\tau\rightarrow 0^+}\widetilde{J}_1^{(2)}(\tau)=0.
\end{align}
For the term $\widetilde{J}_1^{(3)}(\tau)$, we rewrite
\begin{align*}
\widetilde{J}_1^{(3)}(\tau)&=\int_1^\infty \cos\bigl(A_2u+A_3u^\alpha\bigr)\,Q(u,\tau)\,
\mathbf{1}_{[1,M_\tau]}(u)du.
\end{align*}
Since $M_\tau\to\infty$ as $\tau\rightarrow 0$, for every fixed $u\ge1$, $\mathbf1_{[1,M_\tau]}(u)\to1$ as $\tau\rightarrow 0$. It is easy to see that
\begin{align*}
\cos\left(A_2(\tau)u+A_3(\tau)u^\alpha\right)Q(u,\tau)\mathbf1_{[1,M_\tau]}(u)\to
(\alpha-p-1)(\alpha-p-2)u^{\alpha-p-3}\cos u,
\end{align*}
as $\tau\rightarrow 0$. Furthermore, since
\begin{align*}
\int_1^\infty Ku^{\alpha-p-3}\,du<\infty,
\end{align*}
Together with (\ref{u4}) and the dominated convergence theorem, this gives
\begin{align}\label{u9}
\lim_{\tau\to0^+}\widetilde{J}_1^{(3)}(\tau)=(\alpha-p-1)(\alpha-p-2)\int_1^\infty u^{\alpha-p-3}\cos u\,du .
\end{align}
It follows from (\ref{u3}) and (\ref{u7})-(\ref{u9}) that
\begin{align}\label{u10}
\lim_{\tau\to0^+}\widetilde{J}_1(\tau)=-\sin1-(\alpha-p-1)\cos 1-(\alpha-1-p)(\alpha-p-2)\int_1^\infty u^{\alpha-p-3}\cos udu.
\end{align}
For the term $\widetilde{J}_2(\tau)$,
\begin{align}\label{u11}
|\widetilde{J}_2(\tau)|
&\leq \bigl(A_1(\tau)\bigr)^{-\frac{\alpha-p}{\alpha}}
   \int_{(A_1(\tau))^{\frac{1}{\alpha}}M_\tau}^\infty v^{\alpha-p-1} e^{-v^\alpha} dv \nonumber \\
&\leq \bigl(A_1(\tau)\bigr)^{-\frac{\alpha-p}{\alpha}}
   e^{-\frac{A_1(\tau) M_\tau^\alpha}{2}}
   \int_{(A_1(\tau))^{\frac{1}{\alpha}}M_\tau}^\infty v^{\alpha-p-1} e^{-\frac{v^\alpha}{2}} dv \nonumber \\
&\leq \bigl(A_1(\tau)\bigr)^{-\frac{\alpha-p}{\alpha}}
   e^{-\frac{A_1(\tau) M_\tau^\alpha}{2}}
   \int_{0}^\infty v^{\alpha-p-1} e^{-\frac{v^\alpha}{2}} dv\nonumber\\
   &=\bigl(A_1(\tau)\bigr)^{-\frac{\alpha-p}{\alpha}}
   e^{-\frac{A_1(\tau) M_\tau^\alpha}{2}}\frac{1}{\alpha}2^{\frac{\alpha-p}{\alpha}}\Gamma\left(1-\frac{p}{\alpha}\right),
\end{align}
which tends to zero as $\tau\rightarrow 0^+$. By (\ref{J12}), (\ref{u10}), and (\ref{u11}), we get that (\ref{u6}) holds. Substituting (\ref{J11}), and (\ref{u6}) into (\ref{J1}) gives that
\begin{align}\label{u12}
\lim_{\tau\to0^+}J_1(\tau)&=\int_0^1u^{\alpha-p-1}\cos udu-\sin1\nonumber\\
&\quad-(\alpha-p-1)\cos 1-(\alpha-1-p)(\alpha-p-2)\int_1^\infty u^{\alpha-p-3}\cos udu.
\end{align}
Similar to the discussion of (\ref{u12}), we have
\begin{align}
\lim_{\tau\to0^+}J_2(\tau) &= \int_0^1 u^{\alpha-p-1}\sin u\,du + \cos 1 - (\alpha-p-1)\sin 1 \nonumber\\
&\qquad - (\alpha-p-1)(\alpha-p-2)\int_1^\infty u^{\alpha-p-3}\sin udu,  \label{u13}\\
\lim_{\tau\to0^+}J_3(\tau) &= \int_0^1 u^{-p}\sin u\,du + \cos 1 + p\sin 1 - p(p+1)\int_1^\infty u^{-p-2}\sin udu. \label{u14}
\end{align}
Lemma \ref{lem:2.3} implies that
\begin{align}\label{u15}
\lim_{\tau\to0^+}J_1(\tau)=\Gamma(\alpha-p)\cos\frac{(\alpha-p)\pi}{2},\quad \lim_{\tau\to0^+}J_2(\tau)=\Gamma(\alpha-p)\sin\frac{(\alpha-p)\pi}{2}.
\end{align}
Integrating by parts yields that
\begin{align}\label{u16}
\int_1^\infty u^{-p}\sin udu
&= \cos 1 - p\int_1^\infty u^{-p-1}\cos udu \nonumber\\
&= \cos 1 + p\sin 1 - p(p+1)\int_1^\infty u^{-p-2}\sin udu .
\end{align}
Together  with  (\ref{u14}) and the fact that
\begin{align*}
\int_0^\infty u^{-p}\sin udu=\Gamma(1-p)\cos \frac{p\pi}{2},\quad 0<p<\alpha-1,
\end{align*}
this gives
\begin{align}\label{u17}
\lim_{\tau\to0^+}J_3(\tau)=\Gamma(1-p)\cos\frac{p\pi}{2}.
\end{align}
By (\ref{u0}), (\ref{u15}), and (\ref{u17}), we find
\begin{align}\label{u18}
\lim\limits_{\tau\rightarrow 0^+}g'(\tau)&=\frac{2\Gamma(1+p)\sin\frac{\pi p}{2}}{\pi}\left(|b\sigma|^\alpha\Gamma(\alpha-p)\cos\frac{(\alpha-p)\pi}{2}-\widetilde{\beta}|b\sigma|^\alpha\Gamma(\alpha-p)\sin\frac{(\alpha-p)\pi}{2}\right.\nonumber\\
&\qquad\left.-(\theta+b\mu)\Gamma(1-p)\cos\frac{p\pi}{2}\right).
\end{align}
Using the mean value theorem, for any $\tau>0$,
 there exists $\tau_\xi\in(0,\tau)$ such that $g'(\tau_\xi)=\frac{g(\tau)-g(0)}{\tau}$. Letting $\tau\rightarrow 0^+$ gives $\tau_\xi\rightarrow 0^+$. Combining this with (\ref{u18}), we get that (\ref{g0}) holds. This completes the proof.
\end{proof}

\begin{theorem} \label{thm2.1}
If
\begin{align} \label{t2.4}
|b\sigma|^\alpha\left(1-\beta\operatorname{sgn}(b)\tan^2\frac{\pi\alpha}{2}\right)>(\theta+b\mu) \frac{-1}{\Gamma(\alpha)\cos\frac{\alpha\pi}{2}},
\end{align}
then there exists a positive number $\tau^*$ such that for any initial value $x_0 \in \mathbb{R}$, the solution of (\ref{model}) satisfies
\begin{align*}
\limsup_{t \to \infty} \frac{1}{t} \log(|x(t)|) < 0 \quad \text{a.s.}
\end{align*}
provided $\tau\in(0,\tau^*)$. In practice, we can choose a positive number $p\in(0,\alpha-1)$ such that
\begin{align}\label{t2.6}
|b\sigma|^\alpha\left(1-\beta\operatorname{sgn}(b)\tan\frac{\pi\alpha}{2}\tan\frac{(\alpha-p)\pi}{2}\right)>(\theta+b\mu) H_\alpha(p),
\end{align}
where  $\tau^*$ denotes the smallest positive root of the equation
$g(\tau)=1$, $g(\tau)$ is defined in Lemma \ref{le2.4}, and
\begin{align}\label{Halphap}
H_\alpha(p)= -\frac{\Gamma(1 - p) \cos \frac{p\pi}{2}}{\Gamma(\alpha - p) \cos \frac{(\alpha - p)\pi}{2}}.
\end{align}

\end{theorem}

\begin{proof} Consider the function
\begin{align*}
Q(s)=|b\sigma|^\alpha\left(1-\beta\operatorname{sgn}(b)\tan\frac{\pi\alpha}{2}\tan\frac{(\alpha-s)\pi}{2}\right)-(\theta+b\mu) H_\alpha(s),\quad s\in (0,\alpha-1).
\end{align*}
Both $H_\alpha(s)$ and $\tan\frac{(\alpha-s)\pi}{2}$ are continuous, so $Q(s)$ is continuous. Observe that
\begin{align*}
\lim_{s\to 0}\tan \frac{(\alpha-s)\pi}{2}=\tan\frac{\pi\alpha}{2},\quad  \lim_{s\to 0} H_\alpha(s)=\frac{-1}{\Gamma(\alpha) \cos \frac{\alpha\pi}{2}}.
\end{align*}
Consequently,
\begin{align*}
\lim_{s\to 0}Q(s)=|b\sigma|^\alpha\left(1-\beta\operatorname{sgn}(b)\tan^2\frac{\pi\alpha}{2}\right)-(\theta+b\mu) \frac{-1}{\Gamma(\alpha)\cos\frac{\alpha\pi}{2}}.
\end{align*}
By (\ref{t2.4}), we conclude that there exists a sufficiently small $p\in(0,\alpha-1)$ such that $Q(p)>0$, i.e., (\ref{t2.6}) holds. Lemma \ref{le2.4} together with (\ref{t2.6}) yields $g'(0^+)<0$. Since $g(0)=1$ and $g'(0^+)<0$, there exists $\tau_0>0$ such that
$g(\tau_0)<1$. Moreover, Lemma \ref{le2.4} asserts that $\lim_{\tau\to\infty}g(\tau)=\infty$.
Hence, by the intermediate value theorem, the equation
$g(\tau)=1$ has a positive solution. By the continuity of $g$, we have $g(\tau^*)=1$. Therefore,
$\tau^*$ is the smallest positive root of the equation
$g(\tau)=1$. Consequently,
\begin{align}\label{gtau}
g(\tau)<1,\qquad \tau\in(0,\tau^*).
\end{align}

Let us now fix $\tau\in(0,\tau^*)$ arbitrarily. Let $t_k=k\tau$ for $k=0,1,2,\ldots$ and set $x_k=x(t_k)$. For $t\in[t_k,t_{k+1}]$, by the characteristic function of the $\alpha$-stable process, it follows that $x(t)$ can be regarded as the solution to the following equation:
\begin{align*}
dx(t)&=\theta x(t)dt+b x_k dL(t)=(\theta x(t)+b\mu x_k)dt+b\sigma x_k dZ(t)
\end{align*}
with initial value $x(t_k)=x_k$ at the time $t_k$, where $Z(t)$ is a standard $\alpha$-stable L\'evy motion with $Z(1)\sim S_\alpha(1,\beta,0)$. By It\^o's formula, we have
\begin{align}\label{2.17}
x(t)=x_k\left[e^{\theta(t-t_k)}+\frac{b\mu}{\theta}(e^{\theta(t-t_k)}-1)+b\sigma\int_{t_k}^te^{\theta(t-s)} dZ(s)\right].
\end{align}
In particular,
\begin{align}\label{s2.1}
x_{k+1}=x_k\left(e^{\theta \tau}+\frac{b\mu}{\theta}(e^{\theta\tau}-1)+b\sigma \chi_k\right),
\end{align}
where $\chi_k=\int_{t_k}^{t_{k+1}}e^{\theta(t_{k+1}-s)}dZ(s)$. Observe that $\chi_k$ is an $\alpha$-stable random variable with distribution $S_\alpha((\frac{e^{\alpha \theta\tau}-1}{\theta\alpha})^{1/\alpha},\beta,0)$ independent of $x_k$. Hence, we can write $\chi_k=(\frac{e^{\alpha \theta\tau}-1}{\theta\alpha})^{1/\alpha}Z_k$, where $Z_k\sim S_\alpha(1,\beta,0)$ and is independent of $x_k$. From (\ref{s2.1}), we get that
\begin{align}\label{s2.2}
x_{k+1}=x_k\left(e^{\theta \tau}+\frac{b\mu}{\theta}(e^{\theta\tau}-1)+b\sigma \left(\frac{e^{\alpha \theta\tau}-1}{\theta\alpha}\right)^{1/\alpha}Z_k\right).
\end{align}
Consequently,
\begin{align}\label{2.20}
\mathbb{E}|x_{k+1}|^p=\mathbb{E}|x_k|^p\mathbb{E}\left|e^{\theta \tau}+\frac{b\mu}{\theta}(e^{\theta\tau}-1)+b\sigma \left(\frac{e^{\alpha \theta\tau}-1}{\theta\alpha}\right)^{1/\alpha}Z_k\right|^p.
\end{align}
Using Definition \ref{de2.1}, Lemma \ref{le2.1}, and (\ref{G}), we have
\begin{align*}
\mathbb{E}\left|e^{\theta \tau}+\frac{b\mu}{\theta}(e^{\theta\tau}-1)+b\sigma\left(\frac{e^{\alpha \theta\tau}-1}{\theta\alpha}\right)^{1/\alpha}Z_k\right|^p=g(\tau).
\end{align*}
It follows from (\ref{gtau}) that
\begin{align}\label{2.21}
\mathbb{E}\left|e^{\theta \tau}+\frac{b\mu}{\theta}(e^{\theta\tau}-1)+b\sigma\left(\frac{e^{\alpha \theta\tau}-1}{\theta\alpha}\right)^{1/\alpha}Z_k\right|^p\leq 1,\quad \tau\in (0,\tau^*).
\end{align}
Recalling (\ref{gtau}), we get that there is a unique $\varepsilon>0$ such that $g(\tau)=e^{-\varepsilon \tau}$. Combining (\ref{2.20}) with (\ref{2.21}) yields that
\begin{align*}
\mathbb{E}|x_{k+1}|^p\leq\mathbb{E}|x_k|^pe^{-\varepsilon \tau}.
\end{align*}
Since this holds for all $k>0$, we have
\begin{align}\label{2.22}
\mathbb{E}|x_{k+1}|^p\leq |x_0|^pe^{-\varepsilon (k+1)\tau},\quad \forall k\geq 0.
\end{align}
It follows from (\ref{2.17}) that
\begin{align}\label{2.23}
\mathbb{E} \left[ \sup_{t_k \leq t \leq t_{k+1}} |x(t)|^p \right]
= \mathbb{E}|x_k|^p \mathbb{E} \left[ \sup_{t_k \leq t \leq t_{k+1}} \left| e^{\theta(t-t_k)}+\frac{b\mu}{\theta}(e^{\theta(t-t_k)}-1)+b\sigma\int_{t_k}^te^{\theta(t-s)} dZ(s) \right|^p \right].
\end{align}
Observe that
\begin{align*}
&\mathbb{E} \left[ \sup_{t_k \leq t \leq t_{k+1}} \left| e^{\theta(t-t_k)}+\frac{b\mu}{\theta}(e^{\theta(t-t_k)}-1)+b\sigma\int_{t_k}^te^{\theta(t-s)} dZ(s) \right|^p \right]\\
&\leq  e^{p\theta \tau} +\frac{|b\mu|^p}{\theta^p}\left(e^{\theta\tau}-1\right)^p +|b\sigma|^p e^{p\theta t_{k+1}} \mathbb{E} \left[ \sup_{t_k \leq t \leq t_{k+1}} \left| \int_{t_k}^t e^{-\theta s} dZ(s) \right|^p \right].
\end{align*}
Lemma \ref{le2.2}  yields that
\begin{align*}
\mathbb{E} \left[ \sup_{t_k \leq t \leq t_{k+1}} \left| \int_{t_k}^t e^{-\theta s} dZ(s) \right|^p \right]\leq c_2\left(\int_{t_k}^{t_{k+1}}e^{-\alpha\theta s}ds\right)^{\frac{p}{\alpha}}\leq c_2\tau^{\frac{p}{\alpha}}e^{-p\theta t_k}.
\end{align*}
Then,
\begin{align}\label{2.24}
\mathbb{E} \left[ \sup_{t_k \leq t \leq t_{k+1}} \left| e^{\theta(t-t_k)}+\frac{b\mu}{\theta}(e^{\theta(t-t_k)}-1)+b\sigma\int_{t_k}^te^{\theta(t-s)} d\widetilde{Z}(s) \right|^p \right]\leq e^{p\theta\tau}\left(1+\frac{|b\mu|^p}{\theta^p}+c_2|b\sigma|^p\tau^{\frac{p}{\alpha}}\right).
\end{align}
By combining (\ref{2.22})-(\ref{2.24}), we can conclude that
\begin{align*}
\mathbb{E} \left[ \sup_{t_k \leq t \leq t_{k+1}} |x(t)|^p \right]
\leq e^{p\theta\tau}\left(1+\frac{|b\mu|^p}{\theta^p}+c_2|b\sigma|^p\tau^{\frac{p}{\alpha}}\right)|x_0|^pe^{-\varepsilon k\tau},\quad \forall k\geq 0.
\end{align*}
It then follows that
\begin{align*}
\mathbb{P} \left( \sup_{t_k \leq t \leq t_{k+1}} |x(t)|^p \geq e^{-0.5\varepsilon k\tau} \right) &\leq \frac{\mathbb{E} \left[ \sup_{t_k \leq t \leq t_{k+1}} |x(t)|^p \right]}{e^{-0.5\varepsilon k\tau}}\\
&\leq e^{p\theta\tau}\left(1+\frac{|b\mu|^p}{\theta^p}+c_2|b\sigma|^p\tau^{\frac{p}{\alpha}}\right)|x_0|^p e^{-0.5\varepsilon k\tau}, \quad \forall k \geq 0.
\end{align*}
By Borel-Cantelli lemma, it follows that
\begin{align*}
\sup_{t_k \leq t \leq t_{k+1}} |x(t)|^p<e^{-0.5\varepsilon k\tau}
\end{align*}
holds for all but finitely many $k$. That is, for almost all $\omega\in\Omega$, there exists an integer $k_0 = k_0(\omega)$ such that
\begin{align*}
\sup_{t_k \leq t \leq t_{k+1}} |x(t, \omega)|^p < e^{-0.5\varepsilon k\tau}, \quad \forall k \geq k_0(\omega).
\end{align*}
Therefore, for $ t_k \leq t \leq t_{k+1}$ and $k \geq k_0$
\begin{align*}
\frac{1}{t} \log (|x(t, \omega)|) < -\frac{0.5\varepsilon k\tau}{p(k+1)\tau}.
\end{align*}
Letting $t\to \infty$, we get
\begin{align*}
\limsup_{t \to \infty} \frac{1}{t} \log (|x(t, \omega)|) \leq -\frac{\varepsilon}{2p}
\end{align*}
for almost all $\omega \in \Omega$. This completes the proof.
\end{proof}

\begin{remark}
The following two special cases are immediate consequences of Theorem~\ref{thm2.1}.

(i) Suppose that $L$ is strictly stable, i.e., $\mu=0$. Then condition \eqref{t2.4} reduces to
\begin{align*}
|b\sigma|^\alpha
\Bigl(1-\beta\operatorname{sgn}(b)\tan^2\frac{\pi\alpha}{2}\Bigr)>\theta\frac{-1}{\Gamma(\alpha)\cos\frac{\alpha\pi}{2}}.
\end{align*}
Accordingly, Theorem~\ref{thm2.1} remains valid. In practice, one may choose
$p\in(0,\alpha-1)$ such that
\begin{align*}
|b\sigma|^\alpha\Bigl(1-\beta\operatorname{sgn}(b)\tan\frac{\pi\alpha}{2}
\tan\frac{(\alpha-p)\pi}{2}\Bigr)>\theta H_\alpha(p),
\end{align*}
and let $\tau^*$ denote the smallest positive root of
$g(\tau)=1$.

(ii) Suppose that $L$ is symmetric, i.e., $\beta=0$. Then condition \eqref{t2.4} becomes
\begin{align*}
|b\sigma|^\alpha>
(\theta+b\mu)\,
\frac{-1}{\Gamma(\alpha)\cos\frac{\alpha\pi}{2}}.
\end{align*}
Accordingly, Theorem~\ref{thm2.1} remains valid. In practice, one may choose
$p\in(0,\alpha-1)$ such that
\begin{align*}
|b\sigma|^\alpha>
(\theta+b\mu)H_\alpha(p),
\end{align*}
and let $\tau^*$ denote the smallest positive root of
$g(\tau)=1$.
\end{remark}

When $\alpha=2$, $\beta=0$, $\mu=0$, and $\sigma=1$, the
$\alpha$-stable L\'evy motion reduces to a Gaussian process;
precisely,
\begin{align*}
\{Z(t)\}_{t\ge0}\stackrel{d}{=}\{\sqrt2\,B(t)\}_{t\ge0},
\end{align*}
where $B(t)$ is a standard Brownian motion. Consequently, equation \eqref{model} becomes
\begin{align}\label{Bmodel}
dx(t)=\theta x(t)dt+\sqrt{2}bx(t)\left(\left[\frac{t}{\tau}\right]\tau\right)dB(t).
\end{align}

Substituting these parameter into Theorem~\ref{thm2.1},
condition \eqref{t2.4} reduces to $b^2>\theta$, which immediately yields the following corollary.

\begin{corollary}\label{cor2.1}
If $b^2>\theta$, then there exists a positive number $\tau^*$ such that for any initial value $x_0 \in \mathbb{R}$, the solution of (\ref{Bmodel}) satisfies
\begin{align*}
\limsup_{t \to \infty} \frac{1}{t} \log(|x(t)|) < 0 \quad \text{a.s.} \label{t2.5}
\end{align*}
provided $\tau\in(0,\tau^*)$. Indeed, one can choose $p\in(0,1)$ satisfying
\begin{align*}
b^2>\frac{\theta}{1-p},
\end{align*}
and let $\tau^*$ be the smallest positive root of
\begin{equation}\label{inte}
1=\frac{2\Gamma(1+p)\sin\frac{\pi p}{2}}{\pi}
\int_0^\infty
\frac{1-e^{-b^2u^2\frac{e^{2\theta\tau}-1}{2\theta}}
\cos(e^{\theta\tau}u)}{u^{p+1}}du,
\end{equation}
which is obtained from the equation $g(\tau)=1$ by setting
$\alpha=2$.
\end{corollary}

\begin{remark}
Corollary \ref{cor2.1} recovers the Brownian motion stability criterion established by Mao \cite{mao2016}. Indeed, under the parameter correspondence
\begin{align*}
\alpha \text{ (in Mao \cite{mao2016}) }=\theta,
\qquad
\sigma \text{ (in Mao \cite{mao2016}) }=\sqrt{2}\,b,
\end{align*}
the stability condition in Corollary~\ref{cor2.1} coincides with that in Theorem~2.1 of Mao \cite{mao2016}. The main difference lies in the characterization of the critical delay $\tau^*$: here it is determined by the integral equation \eqref{inte}, whereas in Mao \cite{mao2016} it is characterized as the smallest positive root of $H_1(\tau)+H_2(\tau)=0$.
\end{remark}

\section{Numerical Simulation}
In this section, we present a numerical example illustrating the theoretical result established in Theorem~\ref{thm2.1}. Consider equation~(\ref{model}) with parameters $\theta=0.3$ and $b=2$, namely,
\begin{align}\label{ex}
dx(t)=0.3x(t)dt+2x\left(\left[\frac{t}{\tau}\right]\tau\right)dZ(t)
\end{align}
for $t\geq 0$ with the initial $x(0)=0.5$, where $\{Z(t)\}_{t\ge0}$ is a standard
$\alpha$-stable L\'evy motion with parameters $(\alpha,\beta,\mu,\sigma)=(1.5,0.2,0,1)$. Choosing $p=0.003$, we have
\begin{align*}
|b\sigma|^\alpha\left(1-\beta\operatorname{sgn}(b)\tan\frac{\pi\alpha}{2}\tan\frac{(\alpha-p)\pi}{2}\right)-(\theta+b\mu) H_\alpha(p)=1.7754 >0
\end{align*}
Using MATLAB, we find the smallest positive root to $g(\tau)=1$ on $\tau>0$ is $\tau^*=0.488616$. Hence, by Theorem \ref{thm2.1}, if $\tau=0.02$, then the solution of (\ref{ex}) satisfies
\begin{align*}
\limsup_{t \to \infty} \frac{1}{t} \log(|x(t)|) < 0 \quad \text{a.s.}
\end{align*}

Let $x_k=x(k\tau)$ with $\tau=0.02$ and $k\in\mathbb{N}$. By (\ref{s2.2}), it follows that
\begin{align}\label{exam}
x_{k+1}=x_k\left(e^{0.3 \tau}+2 (\frac{e^{0.009}-1}{0.45})^{2/3}\widetilde{Z}_k\right),
\end{align}
where $\widetilde{Z}_k\sim S_{1.5}(1, 0.2, 0)$. We generate the i.i.d random variables $\widetilde{Z}_k$ according to the algorithm of Janicki and Weron \cite{Janicki}. We first generate a random variable
$V_k$ uniformly distributed on $(-\pi/2,\pi/2)$ and an exponentially distributed random variable $W_k$ with mean 1, where $V_k$ and $W_k$ are independent. Then the corresponding skewed $\alpha$-stable random variable is given by
\begin{align*}
\widetilde{Z}_k=D_{\alpha,\beta}\frac{\sin(\alpha (V_k+C_{\alpha,\beta}))}{(\cos(V_k))^{1/\alpha}}\left(\frac{\cos(V_k-\alpha(V_k+C_{\alpha,\beta}))}{W_k}\right)^{(1-\alpha)/\alpha}
\end{align*}
where $C_{\alpha,\beta}=\frac{\arctan(\beta\tan(\pi\alpha/2))}{1-|1-\alpha|}$, $D_{\alpha,\beta}=(\cos(\arctan(\beta\tan(\pi\alpha/2))))^{-1/\alpha}$, $\alpha=1.5$, and $\beta=0.2$. Figure~1(a) displays a sample path of $x(t)$ with initial value
$x_0=1$. Figure~1(b) displays the sample path of $\frac{\ln|x(t)|}{t}$,
illustrating the exponential stability predicted by Theorem~\ref{thm2.1}.

\begin{figure}[htbp]
\centering
\subfigure[]{
\begin{minipage}[t]{0.4\linewidth}
\centering
\includegraphics[width=2.5in]{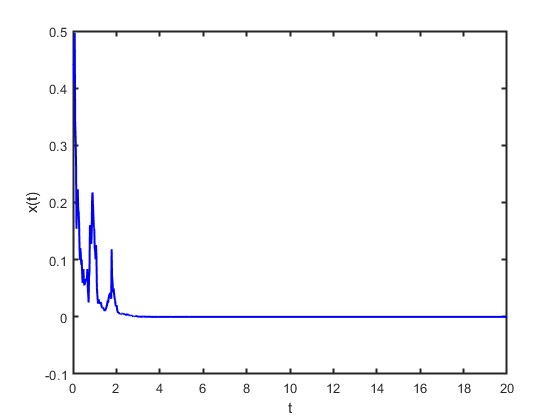}
%\caption{fig1}
\end{minipage}%
}%
%\subfigure[pic2.]{
\subfigure[]{
\begin{minipage}[t]{0.4\linewidth}
\centering
\includegraphics[width=2.5in]{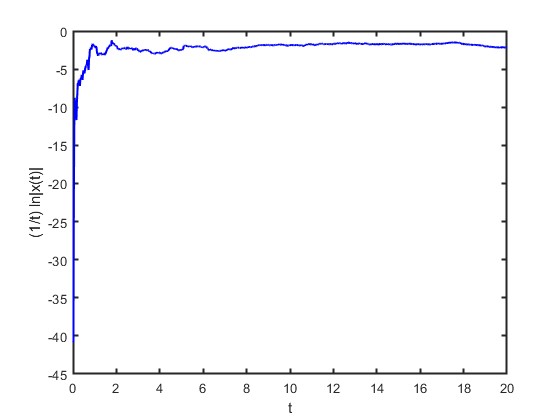}
%\caption{fig2}
\end{minipage}%
}%
\centering\\
Figure 1: \begin{scriptsize}(a) Sample path of $x(t)$ given by (\ref{exam}). (b) Sample path of $\frac{\ln|x(t)|}{t}$  given by (\ref{exam}). \end{scriptsize}
\end{figure}

\section*{Acknowledgments}

The authors are very grateful to the referees and the editor for their insightful and valuable comments which have greatly improved the presentation of the article. This work was supported by the Anhui Provincial Natural Science Foundation (2508085MA016), Anhui Province Action Plan for the Development of Young and Mid-career Faculty in Higher Education Institutions (YQYB2024025), the Natural Science Research Project of Anhui Educational Committee (2023AH030021).

\section*{Data availability}

All data generated or analysed during this study are included in this article.

\section*{Disclosure statement}
No potential conflict of interest was reported by the authors.

\end{document}